\documentclass[11pt,reqno]{article}

\usepackage[usenames]{color}
\usepackage{amssymb}
\usepackage{amsmath}
\usepackage{amsthm}
\usepackage{amsfonts}
\usepackage{amscd}
\usepackage{graphicx}

\usepackage[colorlinks=true,
linkcolor=webgreen,
filecolor=webbrown,
citecolor=webgreen]{hyperref}

\definecolor{webgreen}{rgb}{0,.5,0}
\definecolor{webbrown}{rgb}{.6,0,0}
\usepackage{color}
\usepackage{fullpage}
\usepackage{float}

\usepackage{cases}

\usepackage{graphics,amsmath,amssymb}
\usepackage{amsthm}
\usepackage{amsfonts}
\usepackage{latexsym}

\begin{document}

\theoremstyle{plain}
\newtheorem{theorem}{Theorem}
\newtheorem{remark}{Remark}
\newtheorem{lemma}{Lemma}
\newtheorem{definition}{Definition}
\newtheorem{proposition}{Proposition}
\newtheorem{corollary}{Corollary}

\begin{center}
\vskip 1cm{\Large\bf A refinement of Lang's formula for the sum of powers of integers}
\vskip .2in \large Jos\'{e} Luis Cereceda \\
{\normalsize Collado Villalba, 28400 (Madrid), Spain} \\
\href{mailto:jl.cereceda@movistar.es}{\normalsize{\tt jl.cereceda@movistar.es}}
\end{center}

\begin{abstract}
In 2011, W. Lang derived a novel, explicit formula for the sum of powers of integers $S_k(n) = 1^k + 2^k + \cdots + n^k$ involving simultaneously the Stirling numbers of the first and second kind. In this paper, we first recall and then slightly refine Lang's formula for $S_k(n)$. As it turns out, the refined Lang's formula constitutes a special case of a well-known relationship between the power sums, the elementary symmetric functions, and the complete homogeneous symmetric functions. In addition, we provide several applications of this general relationship.
\end{abstract}

\section{Introduction}

For integers $n \geq 1$ and $k \geq 0$, let $S_k(n)$ denote the sum of the $k$-th powers of the first $n$ positive integers $1^k + 2^k + \cdots + n^k$. In a 2011 technical note \cite{lang}, W.~Lang derived the following explicit formula for $S_k(n)$ (in our notation):
\begin{equation}\label{lang}
S_k(n) = \sum_{m=0}^{\min{(k,n-1)}} (-1)^m (n-m) \genfrac{[}{]}{0pt}{}{n+1}{n+1-m}
\genfrac{\{}{\}}{0pt}{}{n+k-m}{n},
\end{equation}
see \cite[Equation (10)]{lang}, where $\genfrac{[}{]}{0pt}{}{k}{j}$ and $\genfrac{\{}{\}}{0pt}{}{k}{j}$ denote the (unsigned) Stirling numbers of the first and second kind, respectively.

For completeness and for its intrinsic interest, in Section 2, we outline the proof of formula \eqref{lang} as given by Lang. Then, in Section 3, we slightly refine Lang's formula \eqref{lang}. The refinement made essentially amounts to the removal of $n$ from the factor $(n -m)$. In Section 4, we show that the refined Lang's formula arises as a direct consequence of the Newton-Girard identities involving the power sums $S_k(n)$ and the elementary symmetric functions with natural arguments. In Section 5, we point out that, actually, the refined Lang's formula constitutes a special case of a well-known relationship between the power sums, the elementary symmetric functions, and the complete homogeneous symmetric functions. We then look at several applications of this general relationship, which, in the context of this paper, we refer to as the generalized Lang's formula.

\section{Proof of Lang's formula}

Following Lang's own derivation \cite{lang}, next we give a simplified proof sketch of formula \eqref{lang}. We start with the ordinary generating function of $S_k(n)$, i.e.
\begin{equation*}
G_n(x) = \sum_{k=0}^{\infty} (1^k + 2^k + \cdots + n^k) x^k = \sum_{j=1}^n \frac{1}{1 -jx}.
\end{equation*}
This generating function can be rewritten in the form
\begin{equation}\label{proof1}
G_n(x) = \frac{P_n(x)}{\prod_{j=1}^n (1- jx)},
\end{equation}
where $P_n(x)$ is the following polynomial in $x$ of degree $n-1$ with coefficients $P_{n,r}$:
\begin{equation}\label{proof2}
P_n(x) = \sum_{j=1}^n \, \prod_{\substack{l =1 \\ l \neq j}}^n ( 1 -lx) = \sum_{r=0}^{n-1} P_{n,r} x^r.
\end{equation}
Hence, noting that $\frac{1}{\prod_{j=1}^n (1- jx)} = \sum_{m=0}^{\infty} \genfrac{\{}{\}}{0pt}{}{n+m}{n} x^m$, from \eqref{proof1} and \eqref{proof2} it follows that
\begin{equation}\label{proof3}
S_k(n) = \sum_{m=0}^{\min{(k,n-1)}} P_{n,m} \genfrac{\{}{\}}{0pt}{}{n+k-m}{n}.
\end{equation}

Now, as pointed out by Lang \cite{lang}, the elementary symmetric functions $\sigma_m(1,2,\ldots,n)$ enter the scene because we have that
\begin{equation}\label{proof4}
\prod_{j=1}^n (1- jx) = \sum_{m=0}^n (-1)^m \sigma_m(1,2,\ldots,n) x^m,
\end{equation}
with $\sigma_0 =1$. In view of \eqref{proof2} and \eqref{proof4}, it is clear that, by symmetry, $P_n(x)$ must be of the form
\begin{equation*}
P_n(x) = \sum_{m=0}^{n-1} C_{n,m} (-1)^m \sigma_m(1,2,\ldots,n) x^m,
\end{equation*}
for certain positive integer coefficients $C_{n,m}$. Indeed, it can be seen that
\begin{align*}
P_{n,0} & = n, \\
P_{n,1} & = (n-1) (-1) (1 +2+\cdots +n) = (n-1) (-1) \sigma_1(1,2,\ldots,n), \\
P_{n,2} & = (n-2) (1\cdot 2 + 1\cdot 3+ \cdots + (n-1)n) = (n-2) \sigma_2(1,2,\ldots,n), \\[-2mm]
\intertext{and, in general,}
P_{n,m} & = \dfrac{n\binom{n-1}{m}}{\binom{n}{m}} (-1)^m \sigma_m(1,2,\ldots,n) = (n-m) (-1)^m \sigma_m(1,2,\ldots,n),
\end{align*}
so that $C_{n,m} = n-m$, for $m =0,1,\ldots, n-1$.

Therefore, recalling \eqref{proof3}, and invoking the well-known relationship $\sigma_m(1,2,\ldots,n) = \genfrac{[}{]}{0pt}{}{n+1}{n+1-m}$ (see, e.g., \cite[Equation (2.6)]{knuth}), we get \eqref{lang}.

\section{A refinement of Lang's formula}

Having considered Lang's original formula for the sum of powers of integers, we show that this formula can be simplified somewhat. To see this, we write \eqref{lang} in the equivalent form
\begin{align*}
S_k(n) & = n \sum_{m=0}^{\min{(k,n)}} (-1)^m \genfrac{[}{]}{0pt}{}{n+1}{n+1-m}
\genfrac{\{}{\}}{0pt}{}{n+k-m}{n} \\
&  \, \qquad\quad + \sum_{m=1}^{\min{(k,n)}} (-1)^{m-1} \, m \, \genfrac{[}{]}{0pt}{}{n+1}{n+1-m}
\genfrac{\{}{\}}{0pt}{}{n+k-m}{n},
\end{align*}
where the second summation on the right-hand side is zero when $k=0$ or, in other words, it applies for the case that $k \geq 1$. Regarding the first summation, it turns out that
\begin{equation}\label{fsum}
\sum_{m=0}^{\min{(k,n)}} (-1)^m \genfrac{[}{]}{0pt}{}{n+1}{n+1-m}
\genfrac{\{}{\}}{0pt}{}{n+k-m}{n} = \delta_{k,0},
\end{equation}
where $\delta_{k,0}$ is the Kronecker's delta. This is so because
\begin{equation*}
\left( \sum_{i \geq 0} (-1)^i \genfrac{[}{]}{0pt}{}{n+1}{n+1-i} x^i \right)
\left( \sum_{j \geq 0} \genfrac{\{}{\}}{0pt}{}{n+j}{n} x^j \right) =1.
\end{equation*}
Consequently, Lang's original formula \eqref{lang} can be reduced to
\begin{equation}\label{lang2}
S_k(n) = n \, \delta_{k,0} + \sum_{m=1}^{\min{(k,n)}} (-1)^{m-1} \, m \, \genfrac{[}{]}{0pt}{}{n+1}{n+1-m}
\genfrac{\{}{\}}{0pt}{}{n+k-m}{n},
\end{equation}
which holds for any integers $n\geq 1$ and $k\geq 0$, and where, as noted above, the summation on the right-hand side is zero when $k=0$. Moreover, for the general case where $k \geq 1$, formula \eqref{lang2} can in turn be expressed without loss of generality as
\begin{equation}\label{merca}
S_k(n) = \sum_{m=1}^k (-1)^{m-1} \, m \, \genfrac{[}{]}{0pt}{}{n+1}{n+1-m}
\genfrac{\{}{\}}{0pt}{}{n+k-m}{n}, \quad k\geq 1,
\end{equation}
assuming the natural convention that $\genfrac{[}{]}{0pt}{}{n+1}{n+1-m} = \sigma_m(1,2,\ldots,n) =0$ whenever $ m > n$.

\section{Connection with the Newton-Girard identities}

As we shall presently see, the refined Lang's formula for $S_k(n)$ in equation \eqref{merca} can be readily obtained from the Newton-Girard identities (cf.~Exercise 2 of \cite{cere}). Let $\{x_1,x_2,\ldots,x_n \}$ denote a (possibly infinite) set of variables and let $\sigma_m(x_1,x_2,\ldots,x_n)$ denote the corresponding elementary symmetric function. Generally speaking, the Newton-Girard identities are, within the ring of symmetric functions, the connection formulas between the generating sets $\{ \sigma_m(x_1,x_2,\ldots,x_n) \}_{m=1}^k$ and $\{ p_m(x_1,x_2,\ldots,x_n) \}_{m=1}^k$, where $k$ stands for any fixed positive integer and the $p_m$'s stand for the power sums $p_m(x_1,x_2,\ldots,x_n) = x_1^m + x_2^m + \cdots + x_n^m$.

For our purposes here, we focus on the case where $x_i =i$, $\forall i$. Also, to abbreviate the notation, in what follows we write $\sigma_m(1,2,\ldots,n)$ in the shortened form $\sigma_m(n)$. Then, for any given positive integer $m$, the Newton-Girard identities can be formulated as follows (see, e.g., \cite[Equation (5)]{gould} and \cite[Theorem 1.2]{mosse})
\begin{equation}\label{nid}
\sum_{j=1}^{m-1} \overline{\sigma}_{m-j}(n) S_j(n) + S_m(n) + m \overline{\sigma}_{m}(n) =0, \quad m \geq 1,
\end{equation}
where $\overline{\sigma}_{j}(n) = (-1)^j \sigma_j(n)$, and where the summation on the left-hand side is zero when $m=1$. Thus, letting successively $m = 1,2,3,\ldots,k$ in \eqref{nid} gives rise to the following system of $k$ equations in the unknowns $S_1(n), S_2(n), \ldots, S_k(n)$:
\begin{align*}
& S_1(n) = -\overline{\sigma}_{1}(n), \\
& \overline{\sigma}_{1}(n) S_1(n) + S_2(n) = -2 \overline{\sigma}_{2}(n), \\
& \overline{\sigma}_{2}(n) S_1(n) + \overline{\sigma}_{1}(n) S_2(n) + S_3(n) = -3 \overline{\sigma}_{3}(n), \\[-1mm]
& \quad \vdots  \\[-1mm]
& \overline{\sigma}_{k-1}(n) S_1(n) + \overline{\sigma}_{k-2}(n) S_2(n) + \cdots +
\overline{\sigma}_{1}(n) S_{k-1}(n) + S_k(n) = -k \overline{\sigma}_{k}(n),
\end{align*}
which can be expressed in matrix form as
\begin{equation*}
\left( \!\!\begin{array}{ccccc}
1 & 0 & 0 & \cdots & 0 \\
\overline{\sigma}_1(n) & 1 & 0 & \ddots & 0  \\
\overline{\sigma}_2(n) & \overline{\sigma}_1(n) & 1 & \ddots & 0 \\
\vdots & \vdots & \ddots & \ddots & 0  \\[1mm]
\overline{\sigma}_{k-1}(n) & \overline{\sigma}_{k-2}(n) & \cdots & \overline{\sigma}_1(n) & 1
\end{array}\right)
\left( \!\!\begin{array}{c}
S_1(n) \\[1mm]
S_2(n) \\[1mm]
S_3(n) \\[1mm]
\vdots \\[1mm]
S_k(n)
\end{array}\right)
=
\left( \!\!\begin{array}{c}
-\overline{\sigma}_1(n)  \\[1mm]
-2\overline{\sigma}_2(n) \\[1mm]
-3\overline{\sigma}_3(n) \\[1mm]
\vdots \\[1mm]
-k\overline{\sigma}_k(n)
\end{array}\right).
\end{equation*}
On the other hand, it is easily seen that the orthogonality relation in equation \eqref{fsum} is equivalent to the matrix identity
\begin{equation*}
\left( \!\!\begin{array}{ccccc}
1 & 0 & 0 & \cdots & 0 \\
\overline{\sigma}_1(n) & 1 & 0 & \ddots & 0  \\
\overline{\sigma}_2(n) & \overline{\sigma}_1(n) & 1 & \ddots & 0 \\
\vdots & \vdots & \ddots & \ddots & 0  \\
\overline{\sigma}_{k-1}(n) & \overline{\sigma}_{k-2}(n) & \cdots & \overline{\sigma}_1(n) & 1
\end{array} \!\!\right)^{-1}
\! \! =
\left( \!\!\begin{array}{ccccc}
1 & 0 & 0 & \cdots & 0 \\
h_1(n) & 1 & 0 & \ddots & 0  \\
h_2(n) & h_1(n) & 1 & \ddots & 0 \\
\vdots & \vdots & \ddots & \ddots & 0  \\
h_{k-1}(n) & h_{k-2}(n) & \cdots & h_1(n) & 1
\end{array} \!\!\right),
\end{equation*}
where $h_k(n) = \genfrac{\{}{\}}{0pt}{}{n+k}{n}$ and $h_0(n) =1$. Hence, it follows that
\begin{equation*}
\left( \!\!\begin{array}{c}
S_1(n) \\[1mm]
S_2(n) \\[1mm]
S_3(n) \\[1mm]
\vdots \\[1mm]
S_k(n)
\end{array}\right)
=
\left( \!\!\begin{array}{ccccc}
1 & 0 & 0 & \cdots & 0 \\
h_1(n) & 1 & 0 & \ddots & 0  \\
h_2(n) & h_1(n) & 1 & \ddots & 0 \\
\vdots & \vdots & \ddots & \ddots & 0  \\
h_{k-1}(n) & h_{k-2}(n) & \cdots & h_1(n) & 1
\end{array} \!\!\right)
\left( \!\!\begin{array}{c}
-\overline{\sigma}_1(n)  \\[1mm]
-2\overline{\sigma}_2(n) \\[1mm]
-3\overline{\sigma}_3(n) \\[1mm]
\vdots \\[1mm]
-k\overline{\sigma}_k(n)
\end{array}\right).
\end{equation*}
Finally, solving for $S_k(n)$, we obtain
\begin{equation*}
S_k(n) = - \sum_{m=1}^{k} m \overline{\sigma}_m(n) h_{k-m}(n),
\end{equation*}
which is just equation \eqref{merca}.

We conclude this section with the following two remarks.

\begin{remark}
The Newton-Girard identities \eqref{nid} can equally be written as the recurrence formula
\begin{equation*}
S_m(n) = (-1)^{m-1} m \sigma_m(n) - \sum_{j=1}^{m-1} (-1)^j \sigma_j(n) S_{m-j}(n), \quad m \geq 1,
\end{equation*}
giving $S_m(n)$ in terms of $\sigma_1(n),\sigma_2(n),\ldots,\sigma_m(n)$ and the earlier power sums $S_j(n)$, $j =1,2,\ldots,m-1$. This recurrence formula may be compared with the following one appearing in \cite[Remark 3]{cere}:
\begin{equation*}
S_m(n) = m! \binom{n+m}{m+1} - \sum_{j=1}^{m-1} \sigma_j(m-1) S_{m-j}(n), \quad m \geq 1.
\end{equation*}
\end{remark}

\begin{remark}
It is to be noted that the formula for $S_k(n)$ in equation \eqref{merca} was (re)discovered by Merca in \cite[Theorem 1]{merca} by manipulating the formal power series for the Stirling numbers.
\end{remark}

\section{Generalized Lang's formula}

The proof given in the preceding section of formula \eqref{merca} naturally generalizes to arbitrary elementary symmetric functions $\sigma_m(x_1,x_2,\ldots,x_n)$, complete homogenous symmetric functions $h_m(x_1,x_2,\ldots,x_n)$, and associated power sums $p_m(x_1,x_2,\ldots,x_n)$ as follows
\begin{equation}\label{merca2}
p_k(x_1,x_2,\ldots,x_n) = \sum_{m=1}^k (-1)^{m-1} m \sigma_m(x_1,x_2,\ldots,x_n) h_{k-m}(x_1,x_2,\ldots,x_n),
\end{equation}
which becomes equation \eqref{merca} when $x_i = i$, $\forall i$. The generalized Lang's formula \eqref{merca2} constitutes
a well-known result in the theory of symmetric functions (see, e.g., \cite[Proposition 3.2]{egge} and \cite[Lemma 2.1]{merca2} for two recent proofs of formula \eqref{merca2}). Next, we briefly discuss some other applications of it.

Consider first the case in which $x_i =1$, $\forall i$. Then, recalling that $\sigma_m(1,1,\ldots,1) = \binom{n}{m}$ and $h_m(1,1,\ldots,1) = \binom{n+m-1}{m}$, from \eqref{merca2} we obtain the identity
\begin{equation*}
\sum_{m=1}^k (-1)^{m-1} \, m \, \binom{n}{m} \binom{n+k-m-1}{k-m} = n,
\end{equation*}
which holds for any integers $k,n \geq 1$. On the other hand, for integers $1 \leq r \leq n$, it turns out that the $r$-Stirling numbers of the first kind are the elementary symmetric functions of the numbers $r,r+1,\ldots,n$, that is, $\genfrac{[}{]}{0pt}{}{n+1}{n+1-m}_r = \sigma_m(r,r+1,\ldots,n)$; and
the $r$-Stirling numbers of the second kind are the complete symmetric functions of the numbers $r,r+1,\ldots,n$, that is, $\genfrac{\{}{\}}{0pt}{}{n+m}{n}_r = h_m(r,r+1,\ldots,n)$ (see \cite[Section 5]{broder}). Therefore, from \eqref{merca2}, we deduce that
\begin{equation*}
r^k + (r+1)^k + \cdots + n^k = \sum_{m=1}^k (-1)^{m-1} \, m \, \genfrac{[}{]}{0pt}{}{n+1}{n+1-m}_r
\genfrac{\{}{\}}{0pt}{}{n+k-m}{n}_r,
\end{equation*}
where $\genfrac{[}{]}{0pt}{}{n+1}{n+1-m}_r = \sigma_m(r,r+1,\ldots,n) =0$ whenever $m >n+1-r$. Clearly, the last equation reduces to \eqref{merca} when $r=1$. A further generalization of \eqref{merca} in terms of the $r$-Whitney numbers of both kinds and the Bernoulli polynomials can be found in \cite{merca3}.

As another application of equation \eqref{merca2}, we can evaluate the sum of even powers of the first $n$ positive integers by using the fact that (see \cite{merca4})
\begin{equation*}
u(n+1, n+1-m) = (-1)^m \sigma_m(1^2, 2^2, \ldots, n^2),
\end{equation*}
and
\begin{equation*}
U(n+m, n) = h_m(1^2, 2^2, \ldots, n^2),
\end{equation*}
where $u(n,k)$ [respectively, $U(n,k)$] are the central factorial numbers of the first [respectively, second] kind with even indices. Therefore, we have \cite[Theorem 1.1]{merca4}
\begin{equation*}
1^{2k} + 2^{2k} + \cdots + n^{2k} = - \sum_{m=1}^k  m \, u(n+1, n+1-m) U(n+k-m, n).
\end{equation*}
Likewise, noting that (see \cite{merca4})
\begin{equation*}
v(n, n-m) = (-1)^m \sigma_m(1^2, 3^2, \ldots, (2n-1)^2),
\end{equation*}
and
\begin{equation*}
V(n-1+m, n-1) = h_m(1^2, 3^2, \ldots, (2n-1)^2),
\end{equation*}
where $v(n,k)$ [respectively, $V(n,k)$] are the central factorial numbers of the first [respectively, second] kind with odd indices, we can evaluate the sum of even powers of the first $n$ odd integers as follows
\begin{equation*}
1^{2k} + 3^{2k} + \cdots + (2n-1)^{2k} = - \sum_{m=1}^k  m \, v(n, n-m) V(n-1+k-m, n-1).
\end{equation*}

\begin{remark}
Incidentally, the above power sum can alternatively be expressed as the following polynomial in $n$:
\begin{equation*}
1^{2k} + 3^{2k} + \cdots + (2n-1)^{2k} = \frac{2^{2k}}{2k+1} \sum_{j=0}^k \binom{2k+1}{2j+1} B_{2k-2j}
\Big(\frac{1}{2} \Big) \, n^{2j+1},
\end{equation*}
where $B_k(\frac{1}{2})$ denotes the Bernoulli polynomial $B_k(x)$ evaluated at $x = \frac{1}{2}$.
\end{remark}

\begin{table}[bbb]
\centering
\vskip 2pt
\scalebox{1}{
\begin{tabular}{|l|rrrrrrrr|}\hline
$n \backslash  \, j  $ & 0 & 1 & 2 & 3 & 4 & 5 & 6 & 7  \\ \hline
0 & 1 &  &  &  &  &  &  &   \\
1 & 0 & 1 &  &  &  &  &  &   \\
2 & 0 & $-2$ & 1 &  &  &  &  &    \\
3 & 0 & 12 & $-8$ & 1 &  &  &  &   \\
4 & 0 & $-144$ & 108 & $-20$ & 1 &  &  &    \\
5 & 0 & 2880 & $-2304$ & 508 & $-40$ & 1 &  &    \\
6 & 0 & $-86400$ & 72000 & $-17544$ & 1708 & $-70$ & 1 &   \\
7 & 0 & 3628800 & $-3110400$ & 808848 & $-89280$ & 4648 & $-112$ & 1 \\ \hline
\end{tabular}}
\caption{The LS numbers of the first kind, $Ps_n^{(j)}$, up to $n =7$.}\label{tb:1}
\end{table}

\begin{table}[ttt]
\centering
\scalebox{1}{
\begin{tabular}{|l|rrrrrrrr|}\hline
$n \backslash  \, j  $ & 0 & 1 & 2 & 3 & 4 & 5 & 6 & 7  \\ \hline
0 & 1 &  &  &  &  &  &  &   \\
1 & 0 & 1 &  &  &  &  &  &   \\
2 & 0 & 2 & 1 &  &  &  &  &    \\
3 & 0 & 4 & 8 & 1 &  &  &  &   \\
4 & 0 & 8 & 52 & 20 & 1 &  &  &    \\
5 & 0 & 16 & 320 & 292 & 40 & 1 &  &    \\
6 & 0 & 32 & 1936 & 3824 & 1092 & 70 & 1 &   \\
7 & 0 & 64 & 11648 & 47824 & 25664 & 3192 & 112 & 1 \\ \hline
\end{tabular}}
\caption{The LS numbers of the second kind, $PS_n^{(j)}$, up to $n =7$.}\label{tb:2}
\end{table}

Our next application concerns the so-called Legendre-Stirling (LS) numbers of the first and second kind, which, following \cite{andrews}, we denote by $Ps_n^{(j)}$ and $PS_n^{(j)}$, respectively. It is assumed that $n$ and $j$ are non-negative integers fulfilling $0 \leq j \leq n$. Table \ref{tb:1} (\ref{tb:2}) displays the first few LS numbers of the first (second) kind. The LS numbers of the first kind are the elementary symmetric functions of the numbers $2,6,\ldots, n(n+1)$, i.e.
\begin{equation*}
Ps_{n+1}^{(n+1-k)} = (-1)^k \sigma_k(2,6,\ldots, n(n+1)),
\end{equation*}
whereas the LS numbers of the second kind are the complete symmetric functions of the numbers $2,6,\ldots, n(n+1)$, i.e.
\begin{equation*}
PS_{n+k}^{(n)} = h_k(2,6,\ldots, n(n+1)).
\end{equation*}
Equivalently, the above two expressions can be written as
\begin{align*}
Ps_{n+1}^{(n+1-k)} & = (-1)^k 2^k \sigma_k(T_1,T_2,\ldots, T_n), \\[-2mm]
\intertext{and}
PS_{n+k}^{(n)} & = 2^k h_k(T_1,T_2,\ldots, T_n),
\end{align*}
respectively, where $T_n = \frac{1}{2} n(n+1)$ is the $n$-th triangular number. Therefore, we conclude from \eqref{merca2} that
\begin{equation}\label{last1}
T_1^k + T_2^k + \cdots + T_n^k = -\frac{1}{2^k} \sum_{m=1}^k  m \, Ps_{n+1}^{(n+1-m)} PS_{n+k-m}^{(n)}.
\end{equation}
In particular, for $k=1$, we have
\begin{equation*}
T_1 + T_2 + \cdots + T_n = \binom{n+2}{3} = -\frac{1}{2} Ps_{n+1}^{(n)}.
\end{equation*}
Let us observe that the sum of the $k$-th powers of the first $n$ triangular numbers can also be expressed as
\begin{align}\label{last2}
T_1^k + T_2^k + \cdots + T_n^k & = \frac{1}{2^k}  \sum_{j=0}^k \binom{k}{j} S_{k+j}(n) \notag \\
& = \frac{1}{2^k}  \sum_{j=0}^k \binom{k}{j} \frac{B_{k+j+1}(n+1) - B_{k+j+1}(1)}{k+j+1}.
\end{align}
Moreover, Merca showed that, see \cite[Corollary 1.1]{merca5} (in our notation)
\begin{equation}\label{last3}
- \sum_{m=1}^k  m \, Ps_{n+1}^{(n+1-m)} PS_{n+k-m}^{(n)} = \frac{(-1)^k}{(k+1)\binom{2k+2}{k+1}}
+ \sum_{j=0}^k \binom{k}{j} \frac{B_{k+j+1}(n+1)}{k+j+1}.
\end{equation}
Hence, combining \eqref{last1} and \eqref{last3}, and taking into account \eqref{last2}, we obtain the identity
\begin{equation*}
\sum_{j=0}^k (-1)^j \binom{k}{j} \frac{B_{k+j+1}}{k+j+1} =  \frac{1}{(k+1)\binom{2k+2}{k+1}}, \quad k \geq 1,
\end{equation*}
where the $B_k$ are the Bernoulli numbers.

Our last application of the generalized Lang's formula \eqref{merca2} involves the Riemann zeta function at positive even integer arguments, $\zeta(2k)$. According to \cite[Equations (2.1) and (3.1)]{merca6}, we have
\begin{align}
\sigma_k \Big( \frac{1}{1^2}, \frac{1}{2^2}, \frac{1}{3^2}, \ldots \,\Big) & =
\frac{\pi^{2k}}{(2k+1)!}, \label{zeta1} \\[-2mm]
\intertext{and}
h_k \Big( \frac{1}{1^2}, \frac{1}{2^2}, \frac{1}{3^2}, \ldots \,\Big) & = \frac{2^{2k} -2}{2^{2k-1}}
\zeta(2k), \quad k \geq 1. \label{zeta2}
\end{align}
Thus, noting that $\zeta(2k) = \sum_{n=1}^{\infty} \frac{1}{n^{2k}} = p_k\big( \frac{1}{1^2}, \frac{1}{2^2}, \frac{1}{3^2}, \ldots \,\big)$, the substitution of relations \eqref{zeta1} and \eqref{zeta2} into equation \eqref{merca2} yields the recursive formula
\begin{equation}\label{zeta3}
\zeta(2k) = (-1)^{k-1} \frac{k \!\cdot \! \pi^{2k}}{(2k+1)!} + \sum_{m=1}^{k-1} (-1)^{m-1} \frac{2m \cdot \pi^{2m}}{(2m+1)!} \big( 1 - 2^{2(m-k) +1} \big)
\zeta(2k -2m),
\end{equation}
with $k \geq 1$, and where the summation on the right-hand side is zero when $k=1$.

We end this section with the following observations regarding equation \eqref{zeta3}.

\begin{remark}
The recursive formula \eqref{zeta3} was obtained by Merca in \cite[Corollary 4.2]{merca6} by considering the formal power series associated with the symmetric functions \eqref{zeta1} and \eqref{zeta2}. See also \cite{merca7} for a systematic derivation of linear recurrence relations for $\zeta(2k)$.
\end{remark}

\begin{remark}
Since $\sigma_1(x_1, x_2, \ldots\, ) = h_1(x_1, x_2, \ldots \,)$, from \eqref{zeta1} and \eqref{zeta2} we readily obtain that $\zeta(2) = \frac{\pi^2}{6}$. Of course, this result also follows by setting $k =1$ in equation \eqref{zeta3}.
\end{remark}

\begin{remark}
The recursive formula \eqref{zeta3} can be correspondingly expressed in terms of the Bernoulli numbers as follows (cf. \cite[Corollary 5.2]{merca6})
\begin{equation*}
B_{2k} = \frac{2}{2k+1} \sum_{j=1}^k j \binom{2k+1}{2j+1} \bigg( \frac{1}{2^{2k-1}} - \frac{1}{2^{2j}} \bigg)
B_{2k-2j}, \quad k \geq 1.
\end{equation*}
\end{remark}

\section{Conclusion}

In this paper, we have brought to light an outstanding (though largely unnoticed) contribution of W.~Lang to the subject of the sums of powers of integers, namely, his formula for $S_k(n)$ stated in equation \eqref{lang}. We have shown that Lang's original formula \eqref{lang} can be slightly refined so that the integer variable $n$ can be effectively removed from the factor $(n-m)$, as can be seen by looking at formula \eqref{lang2}. Furthermore, we have shown that the refined Lang's formula for $S_k(n)$ in equation \eqref{merca} follows straightforwardly from the Newton-Girard identities formulated in equation \eqref{nid}. Finally, to broaden the scope of the present paper, we have examined several extensions of formula \eqref{merca} achieved by Merca \cite{merca2,merca3,merca4,merca5,merca6}.

Additionally, it should be mentioned that, by considering certain symmetric triples of power series, O'Sullivan \cite{sullivan} provided a natural framework for studying systematically a variety of combinatorial and number theoretic sequences (see, in particular, \cite[Example 5.4]{sullivan}, where the formula in equation \eqref{merca} is obtained by considering a specific symmetric triple).

\end{document}